\documentclass[a4paper,12pt]{amsart}
\usepackage{amsmath,amssymb,latexsym,amsfonts,amscd}

\title{The 5-canonical system on 3-folds of general type}
\author{Jungkai A. Chen, Meng Chen and De-Qi Zhang}
\address{Department of Mathematics, National Taiwan University,
Taipei, 106, Taiwan} \email{jkchen@math.ntu.edu.tw}

\address{School of Mathematical Sciences, Fudan University,
Shanghai, 200433, People's Republic of China  (and Key Laboratory
of Mathematics for Nonlinear Sciences (Fudan University), Ministry
of Education)}
\email{mchen@fudan.edu.cn}

\address{Department of Mathematics, National University of Singapore,
2 Science Drive 2, Singapore 117543, Singapore}
\email{matzdq@nus.edu.sg}

\thanks{The first author was partially supported by the National Science Council and National Center for
Theoretical Science of Taiwan.
 The second author was supported by the National Natural Science
Foundation of China (Key project No.10131010). The third author is
supported by an Academic Research Fund of NUS}

\newcommand{\bP}{{\mathbb P}}
\newcommand{\roundup}[1]{\lceil{#1}\rceil}

\newcommand\CC{{\mathbb{C}}}

\newtheorem{thm}{Theorem}[section]
\newtheorem{lem}[thm]{Lemma}

\newtheorem{prop}[thm]{Proposition}
\newtheorem{claim}[thm]{Claim}
\theoremstyle{definition}

\newtheorem{setup}[thm]{}

\newtheorem{exmp}[thm]{Example}

\newtheorem{rem}[thm]{Remark}
\theoremstyle{remark}

\begin{document}
\begin{abstract}
Let $X$ be a projective minimal Gorenstein 3-fold of general type
with canonical singularities. We prove that the 5-canonical map is
birational onto its image.
\end{abstract}
\maketitle
\pagestyle{myheadings} \markboth{\hfill J. Chen, M. Chen and D.Q.
Zhang\hfill}{\hfill Pluricanonical maps\hfill}
\section{\bf Introduction}
One main goal of algebraic geometry is to classify algebraic
varieties. The successful 3-dimensional MMP (see \cite{KMM, K-M}
for example) has been attracting more and more mathematicians to
the study of algebraic 3-folds. In this paper, we restrict our
interest to projective minimal Gorenstein 3-folds $X$ of general
type where there still remain many open problems.

Denote by $K_X$ the canonical divisor and $\Phi_m :=
\Phi_{|mK_X|}$ the m-canonical map. There has been a lot of work
along the line of the canonical classification. For instance, when
$X$ is a smooth 3-fold of general type with  pluri-genus $h^0(X,
kK_X) \ge 2$, in \cite{Ko}, as an application to his research on
higher direct images of dualizing sheaves, Koll\'ar proved that
$\Phi_{m}$, with $m = 11 k + 5$, is birational onto its image.
This result was improved by the second author \cite{Q-div} to
include the cases $m$ with $m \ge 5k+6$; see also \cite{5map},
\cite{IJM} for results when some additional restrictions (like
bigger $p_g(X)$) are imposed.

On the other hand, for 3-folds $X$ of general type with $q(X)>0$,
Koll\'ar \cite{Ko} first proved that $\Phi_{225}$ is birational.
Recently, the first author and Hacon \cite{CH} proved that
$\Phi_{m}$ is birational for $m \ge 7$ by using the Fourier-Mukai
transform. Moreover, Luo \cite{Lu1}, \cite{Lu2} has some results
for 3-folds of general type with $h^2({\mathcal O}_X)>0$.

Now for minimal and smooth projective 3-folds, it has been
established that $\Phi_{m}$ ($m \ge 6$) is a birational morphism
onto its image after 20 years of research, by Wilson \cite{Wilson}
in  1980, Benveniste \cite{Ben} in  1986 ($m \ge 8$), Matsuki
\cite{Matsuki} in  1986 ($m = 7$), the second author \cite{6map}
in  1998 ($m = 6$) and independently by Lee \cite{Lee1},
\cite{Lee2} in  1999-2000 ($m = 6$; and also the base point
freeness of $m$-canonical system for $m \ge 4$).

The aim of this paper is to prove the following:

\begin{thm}\label{main} Let $X$ be a projective minimal Gorenstein
3-fold of general type with canonical singularities. Then the
m-canonical map $\Phi_m$ is a birational morphism onto its image
for all $m\ge 5$.
\end{thm}

This result is unexpected previously.  The difficulty lies in the
case with smaller $p_g(X)$ or $K_X^3$. One reason to account for
this is that the non-birationality of the $4$-canonical system for
surfaces may happen when they have smaller $p_g$ or $K^2$ (see
Bombieri \cite{Bom}), whence a naive induction on the dimension
does not work.

Nevertheless, there is also evidence supporting the birationality
of $\Phi_5$ for Gorenstein minimal 3-folds $X$ of general type.
For instance, one sees that $K_X^3\ge 2$ for minimal and smooth
$X$ (see \ref{even} below). So an analogy of Fujita's conjecture
would predict that $|5K_X|$ gives a birational map. We recall that
Fujita's conjecture (the freeness part) has been proved by Fujita,
Ein-Lazarsfeld \cite{E-L} and Kawamta \cite{Ka3} when $\dim X \le
4$.

\begin{exmp} The numerical bound "5" in Theorem \ref{main} is optimal.
There are plenty of supporting examples. For instance, let
$f:V\longrightarrow B$ be any fibration where $V$ is a smooth
projective 3-fold of general type and $B$ a smooth curve. Assume
that a general fiber of $f$ has a minimal model $S$ with $K_S^2=1$
and $p_g(S)=2$. (For example, take the product.) Then
$\Phi_{|4K_V|}$ is evidently not birational (see \cite{Bom}).
\end{exmp}

\begin{setup}{\bf Reduction to birationality.}\label{reduction} According to
\cite{6map} or \cite{Lee1}, to prove Theorem \ref{main}, we only
need to verify the statement in the case $m = 5$. On the other
hand, the results in \cite{Lee1, Lee2} show that $|mK_X|$ is base
point free for $m\ge 4$. So it is sufficient for us to verify the
birationality of $|5K_X|$ in this paper.
\end{setup}

\begin{setup}{\bf Reduction to factorial models.} According to the work of M. Reid
\cite{Reid83} and Y. Kawamata \cite{Kafact} (Lemma 5.1), there is
a minimal model $Y$ with a birational morphism
$\nu:Y\longrightarrow X$ such that $K_Y=\nu^*(K_X)$ and that $Y$
is factorial with at worst terminal singularities. Thus it is
sufficient for us to prove Theorem \ref{main} for minimal
factorial models.
\end{setup}

\noindent {\bf Acknowledgments.} We are indebted to H\'el\`ene
Esnault, Christopher Hacon, Yujiro Kawamata, Miles Reid, I-Hsun
Tsai, Eckart Viehweg and Chin-Lung Wang for useful conversations
or comments on this subject. We would like to thank the referee
for a very careful reading and valuable suggestions for
the mathematical and linguistic improvement of the paper.

\section{\bf Notation, Formulae and Set up}
We work over the complex numbers field $\CC$. By {\it a minimal
variety} $X$, we mean one with nef $K_X$ and with terminal
singularities (except when we specify the singularity type).

\begin{setup}\label{begin2}
 Let $X$ be a projective minimal Gorenstein 3-fold of general type.
Take a special resolution $\nu:Y\longrightarrow X$ according to
Reid (\cite{Reid83}) such that $c_2(Y)\cdot \bigtriangleup=0$ (see
Lemma 8.3 of \cite{Miyaoka}) for any exceptional divisor
$\bigtriangleup$ of $\nu$. Write $K_Y=\nu^*K_X+E$ where $E$ is
exceptional and is mapped to a finite number of points. Then for
$m \ge 2$, we have (by the vanishing in \cite{KaV}, \cite{V} or
\cite{E-V}):
$$\chi({\mathcal O}_X)=\chi({\mathcal O}_Y)=-\frac{1}{24}K_Y\cdot
c_2(Y)=-\frac{1}{24}\nu^*K_X\cdot c_2(Y).$$
\begin{eqnarray*}
P_m(X)&= &\chi({\mathcal O}_X(mK_X))=\chi({\mathcal
O}_Y(m\nu^*K_X))\\
&=&\frac{1}{12}m(m-1)(2m-1)K_X^3+\frac{m}{12}\nu^*K_X\cdot
c_2(Y)+\chi({\mathcal O}_Y)\\
&=&(2m-1)(\frac{m(m-1)}{12}K_X^3-\chi({\mathcal O}_{X})).
\end{eqnarray*}
The inequality of Miyaoka and Yau (\cite{Miyaoka}, \cite{Yau})
says that $3c_2(Y)-K_Y^2$ is pseudo-effective. This gives
$\nu^*K_X\cdot(3c_2(Y)-K_Y^2)\ge 0$. Noting that $\nu^*K_X\cdot
E^2=0$ under this situation, we get:
$$-72\chi({\mathcal O}_X)-K_X^3\ge 0.$$
In particular, $\chi({\mathcal O}_X)<0$. So one has:
$$q(X)=h^2({\mathcal O}_X)+(1-p_g(X))-\chi({\mathcal O}_X)>0$$
whenever $p_g(X)\le 1$.
\end{setup}

\begin{setup} \label{even} Suppose that $D$ is any divisor on a smooth 3-fold $V$. The Riemann-Roch
theorem gives:
$$\chi({\mathcal O}_{V}(D))=\frac{D^3}{6}-\frac{K_V\cdot D^2}{4}+\frac{D\cdot(K_V^2+c_2)}{12}+\chi({\mathcal O}_{V}).$$
Direct calculation shows that
$$\chi({\mathcal O}_{V}(D))+\chi({\mathcal O}_{V}(-D))=\frac{-K_V\cdot
D^2}{2}+2\chi({\mathcal O}_{V})\in{\mathbb Z}.$$ Therefore,
$K_V\cdot D^2$ is an even integer.

Now let $X$ be a projective minimal Gorenstein 3-fold of general
type. Let $D$ be any Cartier divisor on X. Then $K_X\cdot
D^2=K_Y\cdot (\nu^*D)^2$ is even. In particular, $K_X^3$ is even
and positive.
\end{setup}

\begin{setup}\label{F} Let $V$ be a smooth projective 3-fold and
let $f:V\longrightarrow B$ be a fibration onto a nonsingular curve
$B$. From the spectral sequence:
$$E_2^{p,q}:=H^p(B,R^qf_*\omega_V)\Longrightarrow
E^n:=H^n(V,\omega_V),$$ Serre duality and Corollary 3.2 and
Proposition 7.6 on pages 186 and 36 of \cite{Ko}, one has the
torsion-freeness of the sheaves $R^i f_* \omega_V$ and the
following:
$$h^2({\mathcal O}_V)=h^1(B,f_*\omega_V)+h^0(B,R^1f_*\omega_V),$$
$$q(V):=h^1({\mathcal O}_V)=g(B)+h^1(B,R^1f_*\omega_V).$$
\end{setup}

\begin{setup}\label{setting} For $\mu=1,2$, we set
$$
\Phi=\begin{cases}
 \Phi_{|K_X|} &\  \text{if}\  p_g(X)\ge 2, \\
 \Phi_{|2K_X|} &\ \text{otherwise}.
 \end{cases}$$
Since we always have $P_2(X)\ge 4$, $\Phi$ is a non-trivial
rational map.

First we fix a divisor $D\in |\mu K_X|$. Let $\pi:
X'\longrightarrow X$ be the composition of both a
desingularization of $X$ and a resolution of the indeterminacy of
$\Phi$. We write $|\pi^*(\mu K_X)|=|M'|+E'$. Then we may assume,
following Hironaka, that:

(1) $X'$ is smooth;

(2) the movable part $M'$ of $|\mu K_{X'}|$ is base point free;

(3) the support of $\pi^*(D)$ is of simple normal crossings.

We will fix some notation below. The frequently used ones are $M$,
$Z$, $S$, $\Delta$ and $E_{\pi}$. Denote by $g$ the composition
$\Phi\circ\pi$. So $g: X'\longrightarrow W'\subseteq{\mathbb P}^N$
is a morphism. Let $g: X'\overset{f}\longrightarrow W
\overset{s}\longrightarrow W'$ be the Stein factorization of $g$
so that $W$ is normal and $f$ has connected fibers. We can write:
$$|\mu K_{X'}|=|\pi^*(\mu K_X)|+\mu E_{\pi}=|M'|+Z',$$
where  $Z'$ is the fixed part and $E_{\pi}$ an effective
$\pi$-exceptional divisor.

 On $X$, one may write $\mu K_X\sim M+Z$
where $M$ is a general member of the movable part and $Z$ the
fixed divisor. Let $S \in |M'|$ be the divisor corresponding to
$M$, then
$$\pi^*(M)=S+\triangle=S+\sum_{i=1}^sd_iE_i$$
with $d_i>0$ for all $i$. The above sum runs over all those
exceptional divisors of $\pi$ that lie over the base locus of
$M$. Obviously $E'=\triangle+\pi^*(Z)$. On the other hand, one may
write $E_{\pi}=\sum_{j=1}^te_jE_j$ where the sum runs over {\it all}
exceptional divisors of $\pi$. One has $e_j>0$ for all $1 \le j \le t$
because $X$ is terminal. Evidently, one has $t\ge s$.

Note that $\text{Sing}(X)$ is a finite set (see \cite{K-M},
Corollary 5.18). We may write $E_{\pi}=\triangle'+\triangle''$
where $\triangle'$ (resp. $\triangle''$) lies (resp. does not lie)
over the base locus of $|M|$. So if one only requires such a
modification $\pi$ that satisfies \ref{setting}(1) and
\ref{setting}(2), one surely has
$\text{supp}(\triangle)=\text{supp}(\triangle').$
\end{setup}

Let $d:=\dim \Phi(X)$. And let $L:=\pi^*(K_X)_{|S}$, which is
clearly nef and big. Then we have the following:

\begin{lem} \label{hit}
When $d\ge 2$, $(L^2)^2 \ge (\pi^*K_X)^3(\pi^*(K_X)\cdot S^2)$.
Moreover, $L^2 \ge 2$.
\end{lem}
\begin{proof}
Take a sufficiently large number $m$ such that $|m\pi^*(K_X)|$ is
base point free. Denote by $H$ a general member of this linear
system. Then $H$ must be a smooth projective surface. On $H$, we
have nef divisors $\pi^*(K_X)_{|H}$ and $S_{|H}$. Applying the
Hodge index theorem, one has
$$(\pi^*(K_X)_{|H}\cdot S_{|H})^2\ge (\pi^*(K_X)_{|H})^2(S_{|H})^2.$$
Removing $m$, we get the first inequality. By \ref{even},
$(\pi^*K_X)^3$ is  even, hence $\ge 2$. Together with
$\pi^*(K_X)\cdot S^2 >0$, we have the second inequality.
\end{proof}

We now state a lemma which will be needed in our proof. The result
might be true for all 3-folds with rational singularities. We
present a proof here just hoping to make this note more
self-contained.

\begin{lem}\label{dq} Let $X$ be a normal projective 3-fold with
only canonical singularities. Let $M$ be a Cartier divisor on $X$.
Assume that $|M|$ is a movable pencil and that $|M|$ has base
points. Then $|M|$ is composed with a rational pencil.
\end{lem}
\begin{proof}
Take a birational morphism $\pi:X'\longrightarrow X$ such that
$X'$ is smooth, that the exceptional divisor $E_{\pi}$ is of
simple normal crossing, and that the map $\Phi_{|M|}$ composed
with $\pi$, becomes a morphism from $X'$ to a curve. Take the
Stein factorization of the latter morphism to get an induced
fibration $f:X'\longrightarrow B$ onto a smooth curve $B$. The
lemma asserts that $B$ must be rational.

Clearly, the exceptional divisor $E_{\pi}$ dominates $B$.

\par \vskip 1pc
{\bf Case 1.} $Bs|M|$ contains a curve $\Gamma$.

This is the easier case. Note that $X$ has only finitely many
points at which $K_X$ is non-Cartier or $X$ is non-cDV (see Cor.
5.40 of \cite{K-M}). So we can pick up a very ample divisor $H$ on
$X$ (avoiding these finitely many points) such that $H$ is Du Val
and intersects $\Gamma$ transversally. We may assume that the
strict transform $H'$ on $X'$ is smooth, i.e., $\pi$ is an
embedded resolution of $H \subset X$. Clearly, there is an
$\pi$-exceptional irreducible divisor $E$ which dominates both
$\Gamma$ and $B$. Now for general $H$, both $H'$ and $E \cap H'$
dominate $B$. Since  the curve $E \cap H'$ arises from the
resolution $\pi : H' \rightarrow H$ of the indeterminancy of the
linear system $|M|_{|H}$ (whose image on $X$ is contained in
$\Gamma \cap H$), it is rational. So $B$ is rational.

\par \vskip 1pc
{\bf Case 2.} $Bs|M|$ is a finite set. (The argument below works
even when $X$ is log terminal.)

Take a base point $P$ of $|M|$. Then
$E = \pi^{-1}(P)$ dominates $B$, i.e., $f(E) = B$.
By Kollar's Theorem 7.6 in \cite{Ko2}, there is an
analytic contractible neighborhood $V$ of $P$ such that $U =
\pi^{-1}(V)\subset X'$ is simply connected. Suppose $g(B)>0$. Then
the universal cover $h : W \longrightarrow B$ of $B$ is
either the
affine line $\CC$ or an open disk in $\CC$.
By Proposition 13.5 of \cite{Fulton}, there is a factorization for
the restriction $f|_U : U \longrightarrow B$, say $ f = h \circ
m$, where $m : U \longrightarrow W$ is continuous.
Note that $m(E)$ is a compact subset of $W$,
so $m(E) $ is a single point. In particular, $f(E)$ is a point, a
contradiction.
\end{proof}

\begin{rem} We received the following
comment about Lemma \ref{dq} from the referee to whom we are much
grateful. Shokurov has already proved that if the pair $(X,
\Delta)$ is klt and the MMP holds, then the fibres of the
exceptional locus are always rationally chain connected, which
easily implies Lemma \ref{dq} in the 3-dimensional case. Further,
the authors noticed that Shokurov's result has recently been
extended by Hacon and McKernan to any dimension and without
assuming MMP.
\end{rem}

\section{\bf The case $p_g\ge 2$}
The following proposition is quite useful throughout the paper.
\begin{prop}\label{key}
Let $S$ be a smooth projective surface. Let $C$ be a smooth curve
on $S$, $N' < N$ divisors on $S$ and $\Lambda \subset |N|$ a
subsystem. Suppose that $|N'|_{|C} = |{N'}_{|C}|$, $\deg(N_{|C}) =
1 + \deg({N'}_{|C}) \ge 1 + 2g(C)$. We consider the following
diagram:
$$
\begin{CD}
|N'| @>{\text{res.}}>> |N'_{|C}| \\
@VV{+\text{eff.}}V  @VV{+P_1}V \\
|N| @>{\text{res.}}>> |N_{|C}| \\
@AA{\subset}A @AA{\subset}A \\
\Lambda @>{\text{res.}}>> \Lambda_{|C} \,\, .
\end{CD}
$$

Suppose furthermore that $\Lambda_{|C}$ is free and $\Lambda_{|C}
\supset |{N'}|_{|C} + P_1$. Then
$$
\Lambda_{|C} = |N|_{|C} = |N_{|C}|, \eqno(*)
$$
which is very ample and complete.
\end{prop}
\begin{proof}
By the Riemann-Roch theorem and Serre duality, we have $\dim$
$|N_{|C}| = 1 + \dim |{N'}_{|C}|$. Since there are inclusions
$|{N'}|_{|C} + P_1 \subseteq \Lambda_{|C} \subseteq |N|_{|C}
\subseteq |N_{|C}|$, now the equalities (*) in the statement
follow from dimension counting and the fact that the first
inclusion above is strict by the freeness of $\Lambda_{|C}$.
\end{proof}

\begin{thm}\label{pg3} Let $X$ be a projective minimal factorial 3-fold of general type.
Assume $p_g(X)\ge 2$. Then $\Phi_5$ is birational.
\end{thm}
\begin{proof}
We give the proof according to the value $d:=\dim \Phi(X)$. As in
\ref{setting}, we set $\Phi=\Phi_1$.

\par \vskip 1pc
{\bf Case 1}: $d=3$. Then $p_g(X)\ge 4$. $\Phi_5$ is birational,
thanks to Theorem 3.1(i) in \cite{IJM}.

\par \vskip 1pc
{\bf Case 2}: $d=2$. We consider the linear system
$|K_{X'}+3\pi^*(K_X)+S|$. Since $K_{X'}+3\pi^*(K_X)+S\ge S$ and
according to Tankeev's principle (see Lemma 2 of \cite{T} or 2.1
of \cite{MPCPS2}), it is sufficient to verify the birationality of
$\Phi_{{|K_{X'}+3\pi^*(K_X)+S|}_{|S}}$. Note that we have a
fibration $f:X'\longrightarrow W$ where a general fiber of $f$ is
a smooth curve $C$ of genus $\ge 2$. The vanishing theorem gives:
$$|K_{X'}+3\pi^*(K_X)+S|_{|S}=|K_S+3L|$$
where $L:=\pi^*(K_X)_{|S}$ is a nef and big divisor on $S$.

By Lemma \ref{hit},  $L^2\ge 2$. According to Reider
(\cite{Reider}), $\Phi_{|K_S+3L|}$ is birational and so is
$\Phi_5$.

\par \vskip 1pc
{\bf Case 3}: $d=1$. In this case, we prefer to replace the
notation $W$ by $B$. Let us set $b:=g(B)$.

Suppose first $b>0$. Let us consider the system $|M|$ on $X$. If
$|M|$ has base points, then $b=0$ by \ref{dq}, a contradiction.
Thus we may assume that $|M|$ is base point free. Then under this
situation $\Phi_5$ is birational, which is exactly the statement
of Theorem 3.3 in \cite{IJM}. We sketch the proof here for the
convenience of the reader. We have an induced fibration
$f:X'\longrightarrow B$. Let $F$ be a general fiber of $f$. Since
$g(B)>0$, the Riemann-Roch and Clifford's Theorem imply that
$S\equiv aF$ with $a\ge p_g(X)\ge 2$. Since $|M|$ is base point
free, one always has $\pi^*(K_X)|_F=\sigma^*(K_{F_0})$ (see Claim
3.3 below), where $\sigma: F\rightarrow F_0$ is the smooth blow
down onto the minimal model. Note that
$$\pi^*(K_X)-F-\frac{1}{a}E'\equiv (1-\frac{1}{a})\pi^*(K_X),$$
which is nef and big. Applying Kawamata-Viehweg vanishing, we have
a surjective map
$$H^0(X', K_{X'}+\roundup{4\pi^*(K_X)-\frac{1}{a}E'})\longrightarrow
H^0(F, K_F+\roundup{(4-\frac{1}{a})\pi^*(K_X)}|_F).$$ Also note
that
$$K_F+\roundup{(4-\frac{1}{a})\pi^*(K_X)}|_F\ge K_F+3\sigma^*(K_{F_0})+
\roundup{(1-\frac{1}{a})E'|_S}.$$ If $(K_{F_0}^2, p_g(F))\neq
(1,2)$, then $|K_F+3\sigma^*(K_{F_0})+
\roundup{(1-\frac{1}{a})E'|_F}|$ defines a birational map by
surface theory and so does $\Phi_{|5K_{X'}|}|_F$. Otherwise, since
$E'|_F\equiv \pi^*(K_X)|_F$ is nef and big, we have the same
conclusion according to \cite{IJM}, Proposition 2.1 which is an
interesting application of Kawamata-Viehweg vanishing and is not
hard to follow. On the other hand, pick up two general fibers
$F_1$ and $F_2$. One has $5K_{X'}\ge
K_{X'}+3\pi^*(K_X)+\bigtriangledown+F_1+F_2$ where
$\bigtriangledown$ is numerically trivial. Kawamata-Viehweg
vanishing gives a surjective map
\begin{eqnarray*}
&&H^0(X',K_{X'}+3\pi^*(K_X)+\bigtriangledown+F_1+F_2)\\
&\longrightarrow& H^0(F_1, K_{F_1}+L_1)\oplus H^0(F_2,
K_{F_2}+L_2),
\end{eqnarray*}
where $L_i:=(3\pi^*(K_X)+\bigtriangledown)_{|F_i}$ is nef and big
for $i=1,2$. Further, the two groups on the right hand side are
non-trivial using Riemann-Roch on the surface $F_i$. This means
that $|5K_{X'}|$ can separate two general fibers of $f$.
Therefore, $\Phi_5$ is birational onto its image.

From now on, we suppose $b=0$. Let $F$ be a general fiber of $f$
and denote by $\sigma:F\longrightarrow F_0$ the smooth blow down onto
the minimal model. We take $\pi$ to be the composition $\pi_1\circ
\pi_0$ where $\pi_0$ satisfies \ref{setting}(1) and
\ref{setting}(2) and $\pi_1$ is a further modification such that
$\pi^*(K_X)$ is supported on a normal crossing divisor.

We may write $S\sim aF$ where $a\ge p_g(X)-1$. And we set
$L:=\pi^*(K_X)_{|F}$ instead. The vanishing theorem gives
$$|K_{X'}+3\pi^*(K_X)+S|_{|F}=|K_F+3L|, $$
from which we see that the problem is reduced to the birationality
of $|K_F+3L|$ because $|K_{X'}+3\pi^*(K_X)+S| \supset |S|$ and
$|S|$ evidently separates different fibers of $f$ (as a line
bundle of positive degree on a rational curve is very ample).
Let $\bar{F}:=\pi_*(F)$. We know that $K_X\cdot \bar{F}^2$ is an
even number by \ref{even}.

If $K_X\cdot \bar{F}^2>0$, then we have
$$L^2=\pi^*(K_X)^2\cdot F=K_X^2\cdot \bar{F}\ge K_X\cdot
\bar{F}^2\ge 2.$$ Reider's theorem says that $|K_F+3L|$ gives a
birational map.

We are left with only the case
$K_X\cdot\bar{F}^2=0$. {}First we have:
\begin{claim}\label{uniform} If $K_X\cdot\bar{F}^2=0$, then
${\mathcal O}_F(\pi^*(K_X)_{|F})\cong {\mathcal
O}_F(\sigma^*K_{F_0}).$
\end{claim}

\begin{proof} It is obvious that the claim is true if it holds
for $\pi=\pi_0$. So we may assume $\pi=\pi_0$. Now
$$0=K_X\cdot (a\bar{F})^2=K_X\cdot M^2=\pi^*(K_X)\cdot
\pi^*(M)\cdot S=a\pi^*(K_X)_{|F}\cdot \triangle_{|F},$$ which
means $\pi^*(K_X)_{|F}\cdot {\triangle'}_{|F}=0$. On the other
hand, the definition of $\pi_0$ gives ${\triangle''}_{|F}=0$. Thus
$(E_{\pi})_{|F}\cdot \pi^*(K_X)_{|F}=0$. The Hodge index theorem
on $F$ tells us that ${E_{\pi}}_{|F}$ must be negative definite.

We may write
$$K_F=\pi^*(K_X)_{|F}+G$$
where $G=(E_{\pi})_{|F}$ is an effective negative definite divisor
on $F$. Note that $L$ is nef and big and that $L\cdot G=0$. The
uniqueness of the Zariski decomposition shows that
$\sigma^*K_{F_0}\sim \pi^*(K_X)_{|F}$. We are done.
\end{proof}

From the above claim, we have $\Phi_{|K_F+3L|}=\Phi_{|4K_F|}$.
We are left to verify the birationality of $\Phi_5$ only when
$\Phi_{|4K_F|}$ fails to be birational, i.e. when $K_{F_0}^2=1$
and $p_g(F)=2$.

Kawamata-Viehweg vanishing  (\cite{E-V, KaV, V}) gives
$$|K_{X'}+3\pi^*(K_X)+F|_{|F}=|K_F+3\sigma^*(K_{F_0})|. \eqno (1)$$
Denote by $C$ a general member of the movable part of
$|\sigma^*K_{F_0}|$. By \cite{BPV}, we know that $C$ is a smooth
curve of genus 2 and $\sigma(C)$ is a general member of
$|K_{F_0}|$. Applying the vanishing theorem again, we have
$$|K_F+2\sigma^*(K_{F_0})+C|_{|C}=|K_{C}+2\sigma^*(K_{F_0})_{|C}|. \eqno (2)$$

Now we may apply Proposition \ref{key}.  Let $N'$ be a divisor
corresponding to the movable part of $|K_F+2\sigma^*(K_{F_0})+C|$
and $N:= (5\pi^*K_X)_{|F}$. Set $\Lambda = |5\pi^*(K_X)|_{|F}$.
It's clear that $N' \le  N$. Also note that $\Lambda$ is free
because $|5K_X|$ is free by \cite{Lee1}.

By $(1)$ above, we see that $\Lambda \supset |N'| +$ (a fixed
effective divisor).

Now restricting to $C$, direct computation shows that $\deg({N'}_{|C})=4$
(by (2)) and $5=\deg(N_{|C}) = 1 + \deg({N'}_{|C})$. Therefore, the induced
inclusion  $|{N'}_{|C}| \hookrightarrow |N_{|C}|$ is given by
adding a single point $P_1$.

By $(2)$, we have $|{N'}_{|C}|=|{N'}|_{|C}$. Together with $(1)$,
we have $\Lambda_{|C} \supset |{N'}_{|C}|+P_1$. Hence by
Proposition \ref{key}, $\Lambda_{|C} = |N_{|C}|$ gives an
embedding. Since $|5 \pi^*K_X|_{|F} \supset |N'| \supset |C|$
(by (1) above) separates different $C$ (noting that $p_g(F) = 2$
and $|C|$ is a rational pencil), ${\Phi_5}_{|F}$ is birational. It
is clear that $|5\pi^* K_X| \supset |S|$ separates different
fibres $F$. Thus $\Phi_5$ is birational.
\end{proof}

\section{\bf Birationality via bicanonical systems}

In this section, we shall complete the proof of Theorem \ref{main}
by studying the bicanonical system.
We set $\Phi:=\Phi_2$ as stated in \ref{setting}. Denote
$d_2:=\dim \Phi_2(X)$. We organize our proof according to the
value of $d_2$.

In the proofs below, we shall apply Tankeev's principle as in the
proof of Theorem 3.2, Case 2.

\begin{thm}\label{d_2=3} Let $X$ be a projective minimal factorial 3-fold of general type.
Assume $d_2=3$. Then $\Phi_5$ is birational.
\end{thm}
\begin{proof}
Recall that $K_X^3$ is even by 2.2, so it's either $>2$  or $=2$.

\par \vskip 1pc
{\bf Case 1}. The case $K_X^3>2$.

Pick up a general member $S$. Let $R:=S_{|S}$. Then $|R|$ is not
composed of a pencil. Thus one obviously has $R^2\ge 2$. So the
Hodge index theorem on $S$ yields
$$\pi^*(K_X)\cdot S^2=\pi^*(K_X)_{|S}\cdot R\ge 2.$$
Set $L:=\pi^*(K_X)_{|S}$. If $K_X^3>2$, then the proof of Lemma
\ref{hit} gives $L^2>2$.

In this case, we must emphasize that we only need  a modification
$\pi$ that satisfies \ref{setting}(1) and \ref{setting}(2).
Namely, we don't need the normal crossings. Thus we have
$\text{Supp}(\bigtriangleup)=\text{Supp}(\triangle').$ This
property is crucial to our proof.

Now the vanishing theorem gives
$$|K_{X'}+2\pi^*(K_X)+S|_{|S} = |K_S+2L|.$$
Since $(2L)^2\ge 12$, we may apply Reider's theorem again.
Assume that $\Phi_{|K_S+2L|}$ is not birational. Then there is a
free pencil $C$ on $S$ such that $L\cdot C=1$. Note that $R\le
2L$, and that $|R|$ is base point free and $|R|$ is not composed
of a pencil. Thus $\dim(\Phi_{|R|}(C))=1$. Since $C$ lies in an
algebraic family and $S$ is of general type, we have $g(C)\ge 2$.
Since $h^0(C, R_{|C})\ge 2$, the Riemann-Roch theorem on $C$ and
Clifford's theorem on $C$ easily imply $R\cdot C\ge 2$. Since
$R\cdot C\le 2L\cdot C=2$, one must have $R\cdot C=2$. Since
$$2L=S_{|S}+\bigtriangleup_{|S}+\pi^*(Z)_{|S}$$ and $C$ is nef, we have
$\bigtriangleup_{|S}\cdot C=0$. This implies that
${\triangle'}_{|S}\cdot C=0$. Note also that
${\triangle''}_{|S}=0$ for  general $S$. We get
$(E_{\pi})_{|S}\cdot C=0$. Therefore
$$K_S\cdot C=(K_{X'}+S)_{|S}\cdot C=\pi^*(K_X)_{|S}\cdot
C+(E_{\pi})_{|S}\cdot C+S_{|S}\cdot C=3,$$ an odd integer. This is
impossible because $C$ is a free pencil on $S$. Therefore, $\Phi_5$ must
be birational.

\par \vskip 1pc
{\bf Case 2}. The case $K_X^3=2$.

If $L^2\ge 3$, then $\phi_5$ is birational according to the proof
in {\bf Case 1}. So we may assume $L^2=2$. By Lemma \ref{hit}, we
have $\pi^*(K_X)\cdot S^2=2$. Set $C = S_{|S}$. Then $|C|$ is base
point free and is not composed with a pencil. So $C^2\ge 2$. The
Hodge index theorem also gives
$$4=(\pi^*(K_X)_{|S}\cdot C)^2\ge L^2\cdot C^2\ge 4.$$
The only possibility is $L^2=C^2=2$ and $L\equiv C$. On the other
hand, the equality
$$4=2K_X^3=K_X^2\cdot (M+Z)=L^2+K_X^2\cdot Z=2+K_X^2\cdot Z$$
gives $K_X^2\cdot Z=2$. Take a very big $m$ such that $|mK_X|$ is
base point free and take a general member $H\in |mK_X|$. By the
Hodge index theorem, $4 = \frac{1}{m^2} (K_X \cdot M \cdot H)^2
\ge \frac{1}{m^2} (K_X^2 \cdot H) (M^2 \cdot H) = 2 K_X \cdot
M^2$. Thus $K_X \cdot M^2 = 2$ and $(K_X)_{|H} \equiv M_{|H}$.
Multiplying by $2$, we deduce that $Z_{|H} \equiv M_{|H}$. Thus
$K_X \cdot Z \cdot M = \frac{1}{m} Z_{|H} \cdot M_{|H} =
\frac{1}{m} M^2 \cdot H = 2$. So $L\cdot \pi^*(Z)_{|S}=2$. Since
$2C \equiv 2L = \pi^*(2K_X)_{|S} = \pi^*(M+Z)_{|S} = (S + \Delta +
\pi^*(Z))_{|S} = C + (\Delta + \pi^*(Z))_{|S}$ and $L^2=L\cdot
C=2$, we see that
$$0 =
L\cdot\Delta = C \cdot \Delta. \eqno(3)$$ Thus $K_S = (K_{X'} +
S)_{|S} = C + (\pi^*(K_X) + E_{\pi})_{|S} = (C + L) +
((E_{\pi})_{|S})=P+N$ is the Zariski decomposition by (3) and
\ref{setting}. Denote by $\sigma: S\longrightarrow S_0$ the
smooth blow down onto the minimal model. Then $C+L\sim
\sigma^*(K_{S_0}).$

Note that $C = S_{|S}$ and $\dim |C| \ge \dim |S|_{|S} \ge 2$
because $|S|$ gives a generically finite map. Assume to the
contrary that $\Phi_5$ is not birational. Then neither is
$\Phi_{|S|}$.
Denote by $d$ the generic degree of $\Phi_5$. Then:
$$2=C^2=S^3\ge d(P_2(X)-3).$$
Because $d\ge 2$, we see $P_2(X)=4$ and $d = 2$.
By the same argument as in Case 1, we have:
$$|5K_{X'}|_{|S}\supset \text{the movable part of }|K_S+2L|\supset |C|,$$
so $\Phi_{|C|}:S\longrightarrow \bP^{h^0(S,C)-1}$ is not
birational either. On the other hand, we may write
$$2=C^2\ge \deg(\Phi_{|C|})\deg(\Phi_{|C|}(S)).$$
If $h^0(S, C)\ge 4$, then $\deg(\Phi_{|C|}(S))\ge 2$ and
$\deg\Phi_{|C|}=1$, i.e. $\Phi_{|C|}$ is birational which
contradicts the assumption. So $h^0(S,C)=3$ and $|C|=|S|_{|S}$.
Therefore, $ \Phi_{|C|}: S\longrightarrow \bP^2$ is generically
finite of degree 2. Let $\Phi_{|C|} = \tau \circ \gamma$ be the Stein
factorization with $\gamma : S \rightarrow S'$ a birational morphism onto
a normal surface and $\tau : S' \rightarrow \bP^2$ a finite morphism of degree 2.
We can write $C = \Phi_{|C|}^*\ell$ with a line $\ell$.

For a curve $E$ on $S$, by the projection formula, $C . E = \ell .
\Phi_{|C|*}E$. So $E$ is contracted to a point on $S'$ if and only
if $E$ is contracted to a point on $\bP^2$ (for $\tau$ is finite);
if and only if $E$ is perpendicular to $C \equiv \frac{1}{2}
\sigma^*(K_{S_0})$ ($=$ half of the pull back of $K_{\overline S}$
which is ample on the unique canonical model $\overline{S}$ of
$S$); if and only if $E$ is contracted to a point on
$\overline{S}$ by the projection formula again; we denote by
$E_{all}$ the union of these $E$. By Zariski's Main Theorem, both
$S \setminus E_{all} \rightarrow \overline{S} \setminus$ (the
image of $E_{all}$) and $S \setminus E_{all} \rightarrow S'
\setminus$ (the image of $E_{all}$) are isomorphisms (so we
identify them). Both $\overline{S}$ and $S'$ are completions of
the same $S \setminus E_{all}$ by adding a finite set. The
normality of $\overline{S}$ and $S'$ implies that the birational
morphisms $S \rightarrow \overline{S}$ and $S \rightarrow S'$  can
be identified, so also $S' = \overline{S}$.

Since $\bar{S}$ is normal, Propositions 5.4, 5.5 and 5.7 of
\cite{K-M} imply a splitting
$$\tau_*{\mathcal O}_{\bar{S}}={\mathcal O}_{\bP^2}\oplus
{\mathcal L}$$ where  ${\mathcal L}$ is a line bundle. Thus we see
that
$$q(S)=q(\bar{S})=h^1(\bar{S},\tau_*{\mathcal
O}_{\bar{S}})=0.$$
Since $S$ is nef and big on $X'$, the long
exact sequence
$$0=H^1(K_{X'}+S)\longrightarrow H^1(K_S)\longrightarrow
H^2(K_{X'})\longrightarrow H^2(K_{X'}+S)=0$$ gives
$q(X)=q(X')=q(S)=0$.  Noting that $\chi({\mathcal O}_X)<0$, we
naturally have $p_g(X)\ge 2$. By Theorem \ref{pg3}, $\Phi_5$ is
birational, a contradiction.

Therefore we have proved the birationality of $\Phi_5$.
\end{proof}

\begin{thm}\label{d_2=2} Let $X$ be a projective minimal factorial 3-fold of general type. Assume $d_2=2$. Then
$\Phi_5$ is birational.
\end{thm}
\begin{proof}
By 2.2, $K_X^3$ is even and hence either $K_X^3 = 2$ or $K_X^3 \ge 4$.

\par \vskip 1pc
{\bf Case 1}. $K_X^3>2$.

When $d_2=2$, $f:X'\longrightarrow W$ is a fibration onto a
surface $W$. Taking a further modification, we may even get a
smooth base $W$. Denote by $C$ a general fiber of $f$. Then
$g(C)\ge 2$. Pick up a general member $S$ which is an irreducible
surface of general type. We may write $S_{|S}\sim
\sum_{i=1}^{a_2}C_i$ where $a_2\ge P_2(X)-2$. Since $K_X^3>2$, we
have $a_2\ge P_2(X)-2\ge 3$. Set $L:=\pi^*(K_X)_{|S}$. Then $L$ is
nef and big. Since $\pi^*(K_X)\cdot S^2=(\pi^*(K_X)_{|S}\cdot
S_{|S})_S\ge 3(\pi^*(K_X)_{|S}\cdot C)_S\ge 3$, Lemma \ref{hit}
gives $L^2\ge 4$. The vanishing theorem gives
$$|K_{X'}+2\pi^*(K_X)+S|_{|S}=|K_S+2L|. \eqno(4) $$

Assume that $\Phi_5$ is not birational. Then neither is
$\Phi_{|K_S+2L|}$ for general $S$. Because $(2L)^2\ge 10$,
Reider's theorem (\cite{Reider}) tells us that there is a free
pencil $C'$ on $S$ such that $L\cdot C'=1$. Since $2 = C'  \cdot
2L \ge C' . S_{|S} = a_2 C' \cdot C \ge 3 C' . C$, we have $C\cdot
C'=0$. So $C'$ lies in the same algebraic family as that of $C$.
We may write
$$2L\equiv a_2C+G$$
where $G=(\Delta+\pi^*(Z))_{|S}\ge 0$ and $a_2\ge 3$.  Since
$2L-C-\frac{1}{a_2}G\equiv (2-\frac{2}{a_2})L$ is nef and big,
Kawamata-Viehweg vanishing gives $H^1(S,
K_S+\roundup{2L-C-\frac{1}{a_2}G})=0$. Thus we get a surjection:
$$H^0(S,K_S+\roundup{2L-\frac{1}{a_2}G})\longrightarrow H^0(C,
K_C+D)$$ where $D:=\roundup{2L-\frac{1}{a_2}G}_{|C}$ with
$\deg(D)\ge (2-\frac{2}{a_2})L\cdot C>1$. Note that $|K_S+2L|
\supset |S_{|S}|$  separates different $C$. If $\deg(D) \ge 3$,
then $|K_C +D|$ defines an embedding, and so does $|K_S + 2L|$, a
contradiction.

So suppose $\deg(D) = 2$. We now apply Proposition \ref{key}. Let
$N'$ be the movable part of $K_S+\roundup{2L-\frac{1}{a_2}G}$ and
let $N = \pi^*(5K_X)_{|S}$. Set $\Lambda:=|5\pi^*(K_X)|_{|S}$.
As in the proof of Theorem \ref{pg3}, we have
$\Lambda \supset |N'| +$ (a fixed effective divisor),
$|N'|_{|C} = |K_C + D|$, $N' \le N$ and $\deg(N_{|C}) = 1 +
\deg(N'_{|C}) = 2g(C) + 1 = 5$ by the calculation:
$$4 \le (2g(C) - 2) + 2= N' \cdot C \le N \cdot C  =  5 \pi^*K_X \cdot C = 5.$$
By Proposition \ref{key},
$\Lambda_{|C} = |N_{|C}|$ gives an
embedding. It is clear that $|5\pi^* K_X| \supset |S|$ separates
different $S$, and $|5 \pi^*K_X|_{|S} (\supset$ the movable
part of $|K_S+2L|$) separates different $C$. Thus $\Phi_5$ is
birational. This is again a contradiction.

\par \vskip 1pc
{\bf Case 2}. $K_X^3=2$.

We first consider the case $L^2\ge 3$. On the surface $S$, we are
reduced to study the linear system $|K_S+2L|$. We have
$$2L\sim S_{|S}+G=\sum_{i=1}^{a_2}C_i+G$$
where $a_2\ge h^0(S, S_{|S}) - 1\ge P_2(X) - 2 \ge 2$. Denote by
$C$ a general fiber of $f:X'\longrightarrow W$. If $a_2\ge 3$, the
proof in {\bf Case 1} already works. So we assume $a_2=2$, then
$P_2(X)=4$, and the image of the fibration
$\Phi_{|S_{|S}|}:S\longrightarrow \bP^{2}$ is a quadric curve
which is a rational curve. This means that $|C|$ is composed with
a rational pencil. Assume that $|K_S+2L|$ does not give a
birational map. Then Reider's theorem says that there is a free
pencil $C'$ on $S$ such that $L\cdot C'=1$. We claim that $C'$ are
$C$ are in the same pencil. In fact, otherwise $C'$ is horizontal with
respect to $C$ and $C\cdot C'>0$. Since $C$ is a rational pencil,
$C\cdot C'\ge 2$. Therefore $L\cdot C'\ge 2$, a contradiction. So
$C'$ lies in the same family as that of $C$ and $L\cdot C=1$. Note
that $K_S+2L=(K_{X'}+2\pi^*(K_X))_{|S}+S_{|S}\ge C$. So $|K_S+2L|$
distinguishes different members in $|C|$. The vanishing theorem
gives
$$H^0(S, K_S+\roundup{2L-\frac{1}{2}G})\longrightarrow H^0(C, K_C+Q)$$
where $Q=\roundup{2L-C-\frac{1}{2}G}_{|C}$ is an effective divisor
on $C$. If $|K_C+Q|$ is not birational, neither is $|K_C|$. So $C$
must be a hyper-elliptic curve and $\Phi_{|K_C|} : C \rightarrow {\bold P}^1$
is a double cover; see Iitaka \cite{Ii}, \S 6.5, page 217.
Suppose $\Phi_5$ is not
birational. (*) Then $\Phi_5$ must be a morphism of generic degree 2.
Set $s=\Phi_5:X\longrightarrow W_5 \subset \bP^N$. Then
$5K_X=s^*(H)$ for a very ample divisor $H$ on the image $W_5$. So
$$5=5\pi^*(K_X)\cdot C=2\deg(H|_{s(\pi(C))})=2\deg_{\bP^N}s(\pi(C))$$
which is a contradiction. Thus $\Phi_5$ must be birational under
this situation.

Next we consider the case $L^2=2$. Lemma \ref{hit} says
$2=\pi^*(K_X)\cdot S^2 = a_2L\cdot C$. We see that $a_2=2$ and
$L\cdot C=1$. We still consider the linear system $|K_S+2L|$. As
above, $a_2=2$ implies that $|C|$ is a rational pencil. Since
$K_S+2L\ge C$, we see that $|K_S+2L|$ distinguishes different
members in $|C|$. By the same argument as above, we have
$$|K_S+2L|_{|C}\supset |K_C+Q|\supset |K_C|.$$
If $\Phi_5$ is not birational, then neither is $\Phi_{|K_S+2L|}$.
This means that $C$ must be a hyper-elliptic curve and $\Phi_5$ is
of generic degree 2. Since $|5K_X|$ is base point
free, we also have a contradiction as in the previous case. So
$\Phi_5$ is birational.
\end{proof}

\begin{thm}\label{d_2=1} Let $X$ be a projective minimal factorial 3-fold of general type.
Assume $d_2=1$. Then $\Phi_5$ is birational.
\end{thm}
\begin{proof} When $X$ is smooth, this theorem has been proved in
\cite{5map}. Our result is a generalization of this result.

Taking a modification $\pi$ as in \ref{setting}, we get an induced
fibration $f:X'\longrightarrow W$ and $B := W$ is a smooth curve
of genus $b:=g(B)$. By Lemma 2.1 of \cite{families}, we know that
$0\le b\le 1$. Let $F$ be a general fiber of $f$.

\begin{claim}\label{equal} We have
$${\mathcal O}_F(\pi^*(K_X)_{|F})\cong {\mathcal
O}_F(\sigma^*(K_{F_0}))$$ where $\sigma: F\longrightarrow F_0$ is
the smooth blow down onto the minimal model.
\end{claim}

\begin{proof}
If $b>0$, then the movable part of $|2K_X|$ is already base point
free by Lemma \ref{dq}. The claim is automatically true.

Suppose $b=0$. Set $\bar{F}:=\pi_*F$. We may write (see \ref{setting}):
$$S=\sum_{i=1}^{a_2}F_i$$
where $a_2\ge P_2(X)-1 \ge 3$ and  $F_i$ is a smooth fiber of $f$
for each $i$. Then $2K_X\equiv a_2\bar{F}+Z$. Assume $K_X\cdot
\bar{F}^2>0$. Then we have
\begin{eqnarray*}
2K_X^3&\ge & a_2K_X^2\cdot\bar{F}\ge a_2^2\\
 &\ge& (P_2(X)-1)^2=\frac{1}{4}(K_X^3-6\chi({\mathcal O}_X)-2)^2\\
 &\ge& \frac{1}{4}(K_X^3+4)^2.
\end{eqnarray*}
The above inequality is absurd. Thus $K_X\cdot \bar{F}^2=0$ and
$\pi^*(K_X)_{|F}\cdot \bigtriangleup_{|F}=0$. Now we apply the
same argument as in the proof of Claim \ref{uniform}. So the
claim is true.
\end{proof}

Considering the linear system $|K_{X'}+2\pi^*(K_X)+S| \supset
|S|$, which evidently separates different fibers of $f$, we get a
surjection by the vanishing theorem:
$$|K_{X'}+2\pi^*(K_X)+S|_{|F}=|K_F+2\sigma^*(K_{F_0})|.$$
Since $F$ is a surface of general type, $\Phi_{|3K_F|}$ is
birational except when $(K_{F_0}^2, p_g(F))=(1,2)$, or $(2,3)$. Thus
$\Phi_5$ is birational except when $F$ is of those two types.

From now on, we assume that $F$ is one of the above two types.
Then $q(F)=0$ according to  surface theory.  By \ref{F}, one has
$q(X)=b$ because $R^1f_*\omega_{X'}=0$. Since we may assume
$p_g(X)\le 1$ by Theorem \ref{pg3} and since $\chi({\mathcal O}_X)<0$ and
$b\le 1$, we see that the only possibility is $q(X)=b=1$,
$p_g(X)=1$ and $h^2({\mathcal O}_X)=0$.

Let $D\in |\pi^*(K_X)|$ be the unique effective divisor. Since
$2D\sim 2\pi^*(K_X),$ there is a hyperplane section $H_2^0$ of
$W'$ in $\bP^{P_2(X)-1}$ such that $g^*(H_2^0) \equiv a_2F$ and
$2D=g^*(H_2^0)+Z'.$ Set
$Z':=Z_v+2Z_h,$ where $Z_v$ is the vertical part with respect to
the fibration $f$ and $2Z_h$ the horizontal part. Thus
$$D=\frac{1}{2} (g^*(H_2^0)+Z_v)+Z_h.$$
Noting that $D$ is a integral divisor, for the general fiber $F$,
$(Z_h)_{|F}=D_{|F}\sim \sigma^*(K_{F_0})$ by Claim 4.4.

Considering the ${\mathbb Q}$-divisor
$$K_{X'}+4\pi^*(K_X)-F-\frac{1}{a_2}Z_v-
\frac{2}{a_2}Z_h,$$ set
$$G:=3\pi^*(K_X)+D-\frac{1}{a_2}Z_v-
\frac{2}{a_2}Z_h$$ and
$$D_0:=\roundup{G}=3\pi^*(K_X)+\roundup{(1-\frac{2}{a_2})Z_h}+
\text{vertical divisors}.$$ For the general fiber $F$, our
$G-F\equiv (4-\frac{2}{a_2})\pi^*(K_X)$ is nef and big. Therefore,
by the vanishing theorem, $H^1(X', K_{X'}+D_0-F)=0$.

We then have a surjective map
$$H^0(X', K_{X'}+D_0)\longrightarrow H^0(F,
K_F+3\sigma^*(K_{F_0})+\roundup{(1-\frac{2}{a_2})Z_h}_{|F}).$$ If
$F$ is a surface with $(K^2, p_g)=(2,3)$, then
$\Phi_{|K_F+3\sigma^*(K_{F_0})+\roundup{(1-\frac{2}{a_2})Z_h}_{|F}|}$
is birational on $F$. Otherwise, since
$$\roundup{(1-\frac{2}{a_2})Z_h}_{|F}\ge
\roundup{(1-\frac{2}{a_2})(Z_h)_{|F}}=\roundup{(1-\frac{2}{a_2})D_{|F}},$$
Proposition 2.1 of \cite{IJM} implies that
$\Phi_{|K_F+3\sigma^*(K_{F_0})+\roundup{(1-\frac{2}{a_2})Z_h}_{|F}|}$
is birational. Thus $\Phi_5$ is birational.
\end{proof}

Theorems \ref{d_2=3},  \ref{d_2=2} and  \ref{d_2=1}, together
with 1.4 and 1.5, imply Theorem
\ref{main}.


\end{document}